\documentclass[a4paper,12pt,oneside]{amsart}
\usepackage{geometry}
\usepackage{epstopdf}
\usepackage{mathtools}
\usepackage[cmtip,all]{xy}
\usepackage{wrapfig}

\usepackage{mathrsfs}
\usepackage{amscd}
\usepackage{amsmath}
\usepackage{amssymb}
\usepackage{amsthm}
\usepackage{float}
\usepackage[dvips]{graphicx}
\usepackage{color}
\usepackage{tikz} 
\usepackage{enumerate}

\usetikzlibrary{matrix,arrows}

\usepackage{array}
\usepackage{amscd}
\usepackage[all]{xy}

\DeclareFontFamily{U}{rsfs}{%
\skewchar\font127}
\DeclareFontShape{U}{rsfs}{m}{n}{%
<-6>rsfs5<6-8.5>rsfs7<8.5->rsfs10}{}
\DeclareSymbolFont{rsfs}{U}{rsfs}{m}{n}
\DeclareSymbolFontAlphabet
{\mathrsfs}{rsfs}
\DeclareRobustCommand*\rsfs{%
\@fontswitch\relax\mathrsfs}

\DeclareFontFamily{U}{rsfs}{%
\skewchar\font127}
\DeclareFontShape{U}{rsfs}{m}{n}{%
<-6>rsfs5<6-8.5>rsfs7<8.5->rsfs10}{}
\DeclareSymbolFont{rsfs}{U}{rsfs}{m}{n}
\DeclareSymbolFontAlphabet
{\mathrsfs}{rsfs}
\DeclareRobustCommand*\rsfs{%
\@fontswitch\relax\mathrsfs}

\theoremstyle{plain}
\newtheorem{theorem}{Theorem}
\newtheorem{thm}{Theorem}[section]

\newtheorem{defi}[thm]{Definition}

\newtheorem{prop*}{Proposition}

\newtheorem{prop-defi}[thm]{Proposition-Definition}
\newtheorem{thm-defi}[thm]{Theorem-Definition}
\newtheorem{lem-defi}[thm]{Lemma-Definition}

\newcommand{\C}{\mathbb{C}}
\newcommand{\ZZ}{\mathbb{Z}}
\newcommand{\Z}{\mathcal{Z}}
\newcommand{\I}{\mathcal{I}}

\newcommand{\pic}{\operatorname{Pic}}
\renewcommand{\O}{\mathcal{O}}
\newcommand{\n}{{\boldsymbol{n}}}

\newcommand{\betab}{{\boldsymbol{\beta}}}

\renewcommand{\S}[1]{S^{[#1]}}

\renewcommand{\i}[1]{\I^{[#1]}}

\newcommand{\dR}{\mathbf{R}}
\newcommand{\dL}{\mathbf{L}}
\renewcommand{\hom}{\mathcal{H}om}
\newcommand{\ext}{\mathcal{E}xt}

\newcommand{\tr}{\operatorname{tr}}

\newcommand{\q}{\operatorname{q}}

\newcommand{\mM}{\mathcal{M}}

\newcommand{\mMw}{\mathcal{M}^{\omega_S}}
\newcommand{\cE}{\mathcal{E}}

\newcommand{\ch}{\operatorname{ch}}
\newcommand{\DT}{\operatorname{DT}}
\newcommand{\VW}{\operatorname{VW}}
\newcommand{\vir}{\operatorname{vir}}
\newcommand{\dT}{\mathbf{T}}

\newcommand{\cT}{\mathcal{T}}
\newcommand{\cZ}{\mathsf Z}

\newcommand{\sP}{\mathsf{P}}

\newcommand{\sA}{\mathsf{A}}

\newcommand{\s}{\mathbf{s}}

\newcommand{\T}{\mathcal{T}}

\newcommand{\pro}{\operatorname{prod}}
\newcommand{\red}{\operatorname{red}}

\newcommand{\cP}{\mathcal{P}}
\newcommand{\sN}{\mathsf{N}}
\newcommand{\sG}{\mathsf{G}}

\newcommand{\sK}{\mathsf{K}}
\newcommand{\sT}{\mathsf{T}}
\newcommand{\sE}{\mathsf{E}}

\newcommand{\cL}{\mathcal{L}}

\newcommand{\PP}{\mathbb{P}}

\title[Nested Hilbert schemes and local DT theory]{Hilbert Schemes, Donaldson-Thomas Theory, Vafa-Witten and Seiberg Witten theories}

\author{Artan Sheshmani}
\begin{document}

\maketitle

\begin{abstract}
This article provides the summary of \cite{GSY17a} and \cite{GSY17b} where the authors studied the enumerative geometry of ``\textit{nested Hilbert schemes}" of points and curves on algebraic surfaces and their connections to threefold theories, and in particular relevant Donaldson-Thomas, Vafa-Witten and Seiberg-Witten theories.  

%When the canonical bundle $\omega_S$ is positive, in combination with Mochizuki's formulas, we are able to express certain equivariant Donaldson-Thomas invariants of stable 2-dimensional sheaves on the total space of $\omega_S$ in terms of our invariants of nested Hilbert schemes, Seiberg-Witten invariants of $S$, and the integrals over the products of Hilbert schemes of points on $S$.
\end{abstract}

\setcounter{tocdepth}{3}
\tableofcontents

\section{Introduction} 

In recent years, there has been extensive mathematical progress in enumerative geometry of surfaces deeply related to physical structures, e.g. around Gopakumar-Vafa invariants (GV); Gromov-Witten (GW), Donaldson-Thomas (DT), as well as Pandharipande-Thomas (PT) invariants of surfaces; and their ``motivic lifts". There is also a tremendous energy in the study of mirror symmetry of surfaces from the mathematics side, especially Homological Mirror Symmetry.  On the other hand, physical dualities in Gauge and String theory, such as Montonen-Olive duality and heterotic/Type II duality have also been a rich source of spectacular predictions about enumerative geometry of moduli spaces on surfaces. For instance an extensive research activity carried out during the past years was to prove the modularity properties of GW or DT invariants as suggested by the heterotic/Type II duality. The first prediction of this type, the Yau-Zaslow conjecture~\cite{MR1398633}, was proven by Klemm-Maulik-Pandharipande-Scheidegger~\cite{KMPS}. Further recent developments in this particularly fruitful direction include: Pandharipande-Thomas proof~\cite{KKV2} of the Katz-Klemm-Vafa conjecture~\cite{KKV} for K3 surfaces, Maulik-Pandharipande proof of modularity of GW invariants for K3 fibered threefolds~\cite{a90}, Gholampour-Sheshmani-Toda proof of modularity of PT and Gholampour-Sheshmani proof of modularity of DT invariants of stable sheaves on K3 fibrations \cite{G-S-Toda, GS18} (moreover, the generalizations of the latter in \cite{sheshmani-diacon}), and finally Gholampur-Sheshmani-Thomas  \cite{G-S-Thomas} proof of modular property of the counting of curves on surfaces deforming freely (in a fixed linear system) in ambient CY threefolds. 

$SU(2)$-\textit{Seiberg-Witten (SW) / DT correspondence}. The introduction of other versions of gauge theories in dimensions 6 started numerous exciting developments in physical mathematics and mathematical physics involving enumerative geometry, mirror symmetry, and related physics of 4 dimensional manifolds, realized as a complex two dimensional subvariety which could exhibit deformations inside ambient higher dimensional target varieties, such as CY threefolds. The SW/DT correspondence is one of such platforms to find interesting structural symmetries between theory of surfaces and theory of threefolds: having fixed a spin structure on a complex surface $S$ the SW invariant of $S$ roughly counts with sign the number of points in the parametrizing space (the moduli space) of solutions to SW equations defined on $S$ \cite[Section 1]{MR1339810}. The focus of SW theory in physics is the study of moduli space of vacua in $N=2$, $D=4$ super Yang-Mills theory, and in particular certain dualities of the theory such as electric-magnetic duality (Montonen-Olive duality). In  \cite{GLSY} Gukov-Liu-Sheshmani-Yau  conjectured a relation between SW invariants (associate to $SU(2)$ gauge invariant theory) of a projective surface $S$ and DT invariants of a noncompact threefold $X$ obtained by total space of a line bundle $L$ on $S$. In a later work \cite{GSY17a, GSY17b} , which will be elaborated below, Gholampour-Sheshmani-Yau  proved the conjecture in \cite[Equation 56]{GLSY}. 

\begin{wrapfigure}{r}{0.2\textwidth}
\includegraphics[width=.2\textwidth]{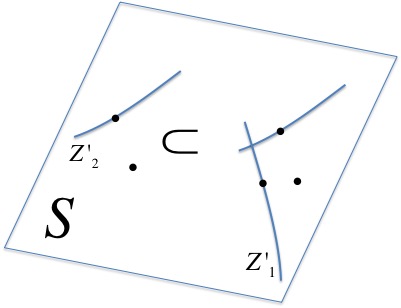}
\end{wrapfigure}
The key object which bridges the geometry of complex surface to the ambient noncompact complex threefold is the ``\textit{nested Hilbert scheme of surface}". This space parametrizes a nested chain of configurations of curves and points in the surface (picture on the right). Let us assume that $S$, as shown in the picture, is a projective simply connected complex surface and let $L$ be a line bundle on $S$. In \cite{GSY17a}, Gholampour-Sheshmani-Yau analyzed the moduli space of stable compactly supported sheaves of modules on the noncompact threefold $X$. In physics terms these are D4-D2-D0 branes wrapping the zero section of $X$ (i.e. sheaves are supported on $S$ or a fat neighborhood of $S$ in $X$). The authors showed in \cite{GSY17a} that the DT invariants of $X$ satisfy an equation in terms of SW invariants of $S$ and certain correction terms governed by invariants of nested Hilbert schemes: 
\begin{center}
DT invariants of $X$= (SW invariants of S)$\cdot$ (Combinatorial coefficients) + Invariants of nested Hilbert scheme of $S$
\end{center}
Their main exciting realization was when $L=K_{S}$, the canonical bundle of $S$. In this case the above DT invariants of $X$ recover some well known invariants in physics, the Vafa-Witten (VW) invariants \cite{VW94}, which are known to have modularity property, following the work of Vafa-Witten \cite{VW94} (and verified by Tanaka-Thomas \cite{TT1,TT2} who defined and computed these invariants mathematically). On the other hand the authors also showed that in some cases the nested Hilbert scheme invariants also have modular property \cite[Theorem 7]{GSY17a}. Therefore, we were able to show that in some cases, the SW invariants can be described in terms of modular forms.

\section*{Acknowledgment}
The author was partially supported by NSF DMS-1607871, NSF DMS-1306313 and Laboratory of Mirror Symmetry NRU HSE, RF Government grant, ag. No 14.641.31.0001. The author would like to further sincerely thank the Center for Mathematical Sciences and Applications at Harvard University, the center for Quantum Geometry of Moduli Spaces at Aarhus University, and the Laboratory of Mirror Symmetry in Higher School of Economics, Russian federation, for the great help and support. 
\section{Nested Hilbert schemes on surfaces}

 Hilbert schemes of points and curves on a nonsingular surface $S$ have been vastly studied. Their rich geometric structures have proved to have many applications in mathematics and physics (see \cite{N99} for a survey).  The current article is a survey of the articles \cite{GSY17a} and \cite{GSY17b}, where the authors studied the enumerative geometry of ``\textit{nested Hilbert schemes}" of points and curves on algebraic surfaces.

Given the sequence  $$\n:=n_{1},n_{2},\dots,n_{r}, r\ge1$$ of nonnegative integers, and $\betab:=\beta_1,\dots,\beta_{r-1}$, a sequence of classes in $H^2(S,\ZZ)$ such that $\beta_i\ge  0$, we denote the corresponding ``\textit{nested}" Hilbert scheme by $\S{\n}_\betab$. A closed point of $\S{\n}_\betab$ corresponds to $$(Z_1,Z_2,\dots,Z_r),\quad (C_1,\dots,C_{r-1})$$ where $Z_i\subset S$ is a 0-dimensional subscheme of length $n_i$, and $C_i\subset S$ is a divisor with $[C_i]=\beta_i$, and  $Z_{i+1}$ is a subscheme of $Z_i\cup C_i$ for any $i<r$, or equivalently \begin{equation}\label{incl} I_{Z_{i}}(-C_i)\subset I_{Z_{i+1}}.\end{equation}

In \cite{GSY17a} and \cite{GSY17b}, in order to define invariants for the nested Hilbert schemes (see \cite[Definitions 2.13, 2.14]{GSY17a}), the authors constructed a virtual fundamental class $[\S{\n}_\betab]^{\vir}$ and then considered cases where one integrates appropriate cohomology classes against it. More precisely, they  constructed a natural perfect obstruction theory over  $\S{\n}_\betab$. This is done by studying the deformation/obstruction theory of the maps of coherent sheaves given by the natural inclusions \eqref{incl} following Illusie.  It turned out that this construction in particular provides a uniform way of studying all known obstruction theories of the Hilbert schemes of points and curves, as well as the stable pair moduli spaces on $S$. The first main result of  \cite{GSY17a} is as follows (see also \cite[Propositions 2.5, and Corollary 2.6]{GSY17a}):
\begin{theorem} \label{thm1}
(\cite[Theorem 1]{GSY17a}) Let $S$ be a nonsingular projective surface over $\C$ and $\omega_S$ be its canonical bundle.The nested Hilbert scheme $\S{\n}_\betab$ with $r\ge 2$ carries a natural virtual fundamental class $$[\S{\n}_\betab]^{\vir} \in A_d( \S{\n}_\betab), \quad \quad  d=n_1+n_r+\frac{1}{2}\sum_{i=1}^{r-1}\beta_i\cdot(\beta_i-c_1(\omega_S)).$$ 
\end{theorem}

\subsection{Special cases}
In the simplest special case, i.e. when $r=1$, we have $\S{\n}_\betab=\S{n_1}$ is the Hilbert scheme of $n_1$ points on $S$ which is nonsingular of dimension $2n_1$, and hence it has a well-defined fundamental class $[\S{n_1}]\in A_{2n_1}(\S{n_1})$. For $r>1$ and $\beta_i=0$, $\S{\n}:=\S{\n}_{(0,\dots,0)}$ is the nested Hilbert scheme of points on $S$ parameterizing the flag of 0-dimensional subschemes $$Z_r\subset \dots \subset Z_2\subset Z_1\subset S,$$ which is in general singular of actual dimension $2n_1$. 

The authors were specifically interested in the case $r=2$, that is: $\S{\n}_\betab=\S{n_1,n_2}_\beta$ for some $\beta\in H^2(S,\ZZ)$.   Interestingly, in the following cases the invariants of nested Hilbert schemes coincide with the Poincar\'e and the stable pair invariants  of $S$ that were previously studied in the context of algebraic Seiberg-Witten invariants and curve counting problems. The following theorem is proven in \cite[Section 3]{GSY17a}.
\begin{theorem}\label{thm1.2} (\cite[Proposition 3.1]{GSY17a}) The virtual fundamental class of Theorem \ref{thm1} recovers the following known cases:
\begin{enumerate}[1.]
\item If $\beta=0$ and $n_1=n_2=n$ then $\S{n,n}_{\beta=0}\cong \S{n}$ and $[\S{n,n}_{\beta=0}]^{\vir}=[\S{n}]$ is the fundamental class of the Hilbert scheme of $n$ points.
\item If $\beta=0$ and $n_2=0$, then $\S{n,0}_{\beta=0}\cong \S{n}$ and $$[\S{n,0}_{\beta=0}]^{\vir}=(-1)^{n}[\S{n}]\cap c_{n}(\omega_S^{[n]}),$$ where $\omega_S^{[n]}$ is the rank $n$ tautological vector bundle over $\S{n}$ associated to the canonical bundle $\omega_S$ of $S$.\footnote{We were notified about this identity by Richard Thomas.}
\item  If $\beta=0$ and $n=n_2=n_1-1$, then it is known that $\S{n+1,n}_{\beta=0}\cong \PP(\i{n})$ is nonsingular, where $\i{n}$ is the universal ideal sheaf over $\S{n}\times S$  \cite[Section 1.2]{L99}. Then, $$[\S{n+1,n}_{\beta=0}]^{\vir}=-[\S{n+1,n}_{\beta=0}]\cap c_1(\O_{\PP}(1)\boxtimes \omega_S).$$
\item If $n_1=n_2=0$ and $\beta\neq 0$, then $\S{0,0}_\beta$ is the Hilbert scheme of divisors in class $\beta$, and $[\S{0,0}_\beta]^{\vir}$ coincides with virtual cycle used to define Poincar\'e invariants in \cite{DKO07}.
\item If $n_1=0$ and $\beta\neq 0$, then $\S{0,n_2}_\beta$ is the relative Hilbert scheme of points on the universal divisor over $\S{0,0}_\beta$, which as shown in \cite{PT10}, is a moduli space of stable pairs and $[\S{0,n_2}_\beta]^{\vir}$ is the same as the virtual fundamental class constructed in \cite{KT14} in the context of stable pair theory. If $P_g(S)=0$ this class was used in  \cite{KT14} to define stable pair invariants.
\end{enumerate}
\end{theorem}

In certain cases, the authors constructed a reduced virtual fundamental class for $\S{n_1,n_2}_\beta$ by reducing the perfect obstruction theory leading to Theorem \ref{thm1}  (\cite[Propositions 2.10, 2.12]{GSY17a}):

\begin{theorem} \label{thm1.5}
(\cite[Theorem 3]{GSY17a}) Let $S$ be a nonsingular projective surface with $p_g(S)>0$, and the class $\beta$ be such that the natural map
$$H^1(T_S)\xrightarrow{*\cup \beta} H^2(\O_S)\quad \text{is surjective,}$$ then, $[\S{n_1,n_2}_\beta]^{\vir}=0$. In this case the nested Hilbert scheme $\S{n_1,n_2}_\beta$ carries a reduced virtual fundamental class $$[\S{n_1,n_2}_\beta]_{\red}^{\vir} \in A_d( \S{n_1,n_2}_\beta), \quad \quad  d=n_1+n_2+\frac{1}{2}\beta\cdot(\beta-K_S)+p_g.$$
The reduced virtual fundamental classes $[\S{0,0}_\beta]^{\vir}_{\red}$ and $[\S{0,n_2}_\beta]^{\vir}_{\red}$   match with the reduced virtual cycles constructed in \cite{DKO07, KT14} in cases 3 and 4 of Proposition \ref{thm1.2}. $[\S{0,n_2}_\beta]^{\vir}_{\red}$ was used in \cite{KT14} to define the stable pair invariants of $S$ in this case.
\end{theorem}

\subsection{Nested Hilbert scheme of points} One interesting specialization of the nested Hilbert schemes is the case where $ \betab=0$ and $\n=n_{1}\geq n_2$. In \cite{GSY17a} the authors studied the nested Hilbert schemes of points $$\S{n_1\ge n_2}:=\S{n_1,n_2}_{\beta=0}$$ in much more details. Let  $\iota: \S{n_1\ge n_2}\hookrightarrow \S{n_1}\times \S{n_2}$ be the natural inclusion. For the case where $S$ is toric with the torus $\dT$ and the fixed set $S^\dT$, the authors provided a purely combinatorial formula for computing $[\S{n_1\ge n_2}]^{\vir}$ by torus localization along the lines of \cite{MNOP06}:

\begin{theorem} \label{thm1.7}
(\cite[Theorem 4]{GSY17a}) For a toric nonsingular surface $S$ the $\dT$-fixed set of $\S{n_1,n_2}$ is isolated and given by tuple of nested partitions of $n_2, n_1$: $$\left\{(\mu_{2,P} \subseteq \mu_{1,P})_P\mid P \in S^\dT,  \quad \mu_{i,P} \vdash n_i \right\}.$$
Moreover, the $\dT$-character of the virtual tangent bundle $\T$ of $\S{n_1\ge n_2}$ at the fixed point $Q=(\mu_{2,P} \subseteq \mu_{1,P})_P$ is given by  $$ \tr_{\cT^{\vir}_{Q}}(t_1,t_2)=\sum_{P\in S^\dT} \mathsf{V}_P,$$ where $t_1, t_2$ are the torus characters and $\mathsf{V}_P$ is a Laurent polynomial in $t_1, t_2$ that is completely determined by the $\mu_{2,P},\mu_{1,P}$ and is given by the right hand side of the following formula \begin{equation} \label{virtan}\tr_{\cT^{\vir}_{I_1\subseteq I_2}}=\cZ_1+\frac{\overline{\cZ}_2}{t_1t_2}+\left(\overline{\cZ}_1\cdot\cZ_2-\overline{\cZ}_1\cdot \cZ_1-\overline{\cZ}_2\cdot \cZ_2\right)\frac{(1-t_{1})(1-t_{2})}{t_{1}t_{2}}.\end{equation}

\end{theorem}

It turns out that when $S$ is toric and Fano, by torus localization, one can express $[\S{n_1\ge n_2}]^{\vir}$ in terms of the fundamental class of the product of Hilbert schemes $\S{n_1}\times \S{n_2}$:

\begin{theorem} \label{thm2} (\cite[Theorem 5]{GSY17a}) If $S$ is a nonsingular projective toric Fano surface, then, $$\iota_*[\S{n_1\ge n_2}]^{\vir}=[\S{n_1}\times \S{n_2}]\cap c_{n_1+n_2}(\sE^{n_1,n_2}),$$ where $\sE^{n_1,n_2}$ is the rank $n_1+n_2$ vector bundle on $\S{n_1}\times \S{n_2}$ obtained by the first relative extension sheaf of the universal ideal sheaves $\i{n_1}$ and $\i{n_2}$.
\end{theorem}

Theorem \ref{thm2} holds in particular for $S=\mathbb{P}^2,\; \mathbb{P}^1\times \mathbb{P}^1$, which are the generators of the cobordism ring of nonsingular projective surfaces. The authors used this fact together with a degeneration formula developed for $[\S{n_1\ge n_2}]^{\vir}$ (\cite[Proposition 5.1]{GSY17a}) to prove:

\begin{theorem} \label{thm3} 
 If $S$ is a nonsingular projective surface, and $\alpha$ is a cohomology class in $H^{n_1+n_2}(\S{n_1}\times \S{n_2})$ with the following properties: 
\begin{itemize}
 \item $\alpha$ is universally defined for any pair of a nonsingular projective surface and a line bundle on it,
 \item $\alpha$ is well-behaved under good degenerations of $S$,
 \end{itemize}
 then $$\int_{[\S{n_1\ge n_2}]^{\vir}}\iota^*\alpha =\int_{\S{n_1}\times \S{n_2}}\alpha \cup c_{n_1+n_2}(\sE^{n_1,n_2}),$$ and $\sE^{n_1,n_2}$ is the rank $n_1+n_2$ virtual vector bundle on $\S{n_1}\times \S{n_2}$ obtained by taking the alternating sum (in the $K$-group) of all the relative extension sheaves of $\i{n_1}$ and $\i{n_2}$.
\end{theorem}

The operators $$\int_{\S{n_1}\times \S{n_2}}- \cup c_{n_1+n_2}(\sE^{n_1,n_2}_M)$$ were studied by Carlsson-Okounkov in \cite{CO12}. Here $M\in \pic(S)$, and $\sE^{n_1,n_2}_M$ is the rank $n_1+n_2$ virtual vector bundle over $\S{n_1}\times \S{n_2}$ obtained by taking the alternating sum of all the relative extensions of $\i{n_1}$ and $\i{n_2}\boxtimes M$ defined as follows:

\begin{defi} \label{virbdl} (\cite[Definition 4.3]{GSY17a}) For any line bundles $M$ on $S$, let $\sE_M^{n_1,n_2} \in K(\S{n_1}\times \S{n_2})$  be the element of virtual rank $n_1+n_2$ defined by $$\sE_M^{n_1,n_2}:=\left[\dR\pi'_*p^*M \right]-\left[\dR\hom_{\pi'}(\i{n_1},\i{n_2}\otimes p^*M)\right], $$
where $p$ and  $\pi'$ are respectively the projections from $S\times \S{n_1}\times \S{n_2}$ to the first and the product of last two factors. Let $i$ be the inclusion of the closed point $(I_1,I_2) \in \S{n_1}\times \S{n_2}$, then, we define
$$\sE^{n_1,n_2}_M|_{(I_1,I_2)}:=\left[\dL i^*\dR\pi'_*p^*M \right]-\left[\dL i^*\dR\hom_{\pi'}(\i{n_1},\i{n_2}\otimes p^*M)\right]\in K(\operatorname{Spec}(\C)).$$
If $M=\O$, we sometimes drop it from the notation. We also define the following generating series  
$$Z_{\pro}(S,M):=\sum_{n_1\ge n_2\ge 0}q_1^{n_1}q_2^{n_2}\int_{\S{n_1}\times \S{n_2}}c(\sE^{n_1,n_2})\cup c(\sE_M^{n_1,n_2}).$$
Note that by equation (5) in \cite{CO12}, we know  $c_i(\sE^{n_1,n_2}_M)=0$ for $i> n_1+n_2$, and so the integrand in the definition of  $Z_{\pro}(S,M)$ can be replaced by $$c_{n_1+n_2}(\sE^{n_1,n_2})\cup c_{n_1+n_2}(\sE_M^{n_1,n_2}).$$
\end{defi}

Carlsson and Okounkov  were able to express these operators in terms of explicit vertex operators. As an application of Theorem \ref{thm3} and the result of \cite{CO12}, the following explicit formula was proven:

\begin{theorem}(\cite[Theorem 7]{GSY17a}) \label{thm4}  Let $S$ be a nonsingular projective surface, $\omega_S$ be its canonical bundle, and $K_S=c_1(\omega_S)$. Then,
\begin{align*}\sum_{n_1\ge n_2\ge 0}(-1)^{n_1+n_2}&\int_{[\S{n_1\ge n_2}]^{\vir}} \iota^*c(\sE^{n_1,n_2}_M)q_{1}^{n_1}q_2^{n_2}=\\&\prod_{n> 0}\left(1- q_2^{n-1}q_1^n\right)^{\langle K_S, K_S-M \rangle}\left(1- q_1^nq_2^{n}\right)^{\langle K_S-M, M \rangle-e(S)},
\end{align*}
where $\langle-,-\rangle$ is the Poincar\'e paring on $S$ and $\sE^{n_1,n_2}_M$ is as in Theorem \ref{thm2}.%, and $\left[q_{1}^{n_1}q_2^{n_2}\right]$ means the coefficient of $q_{1}^{n_1}q_2^{n_2}$ in the expansion of the expression on the right hand side.
\end{theorem}
\section{Nested Hilbert schemes and DT theory of local surfaces}
In this section we give an overview of the results of \cite{GSY17b}. Let $(S,h)$ be a nonsingular simply connected projective surface with $h=c_1(\O_S(1))$. Let $\omega_S$ be the canonical bundle of $S$ with the projection map $\q$ to $S$, and $X$ be the total space of $\omega_S$ \footnote{In \cite{GSY17b}, we consider a more general case in which $X$ is the total space of an arbitrary line bundle $\cL$ with $H^0(\cL\otimes \omega_S^{-1})\neq 0$.}. $X$ is a noncompact Calabi-Yau threefold and one can define the DT invariants of $X$ by using $\C^*$-localization, where $\C^*$-acts on $X$ by scaling the fibers of $\omega_S$. More precisely, let $$v=(r,\gamma,m) \in \oplus_{i=0}^2 H^{2i}(S,\mathbb{Q})$$ be a Chern character vector,  and $\mMw_h(v)$ be the moduli space of compactly supported 2-dimensional stable sheaves $\cE$ on $X$ such that $\ch(\q_*\cE)=v$. Here stability is defined by means of the slope of $\q_*\cE$ with respect to the polarization $h$.  In \cite{GSY17b} the authors provided $\mMw_h(v)$ with a perfect obstruction theory by reducing the natural perfect obstruction theory given by \cite{T98}. The fixed locus $\mMw_h(v)^{\C^*}$ of the moduli space is compact and the reduced obstruction theory gives a virtual fundamental class over it, denoted by $[\mMw_h(v)^{\C^*}]^{\vir}_{\red}$. They defined two types of 
DT invariants:

\begin{align*}
\DT^{\omega_S}_h(v;\alpha)&=\int_{[\mMw_h(v)^{\C^*}]^{\vir}_{\red}} \frac{1}{e(\text{Nor}^{\vir})}\in \mathbb{Q}[\s,\s^{-1}],\quad \quad \alpha \in H^*_{\C^*}(\mMw_h(v)^{\C^*})_\s\\
\DT^{\omega_S}_h(v)&=\chi^{vir}(\mMw_h(v)^{\C^*})\in \ZZ,
\end{align*} where $\text{Nor}^{\vir}$ is the virtual normal bundle of $\mMw_h(v)^{\C^*}\subset \mMw_h(v)$, $\chi^{vir}(-)$ is the virtual Euler characteristic \cite{FG10}, and $\s$ is the equivariant parameter. 

If $\alpha=1$ then the authors were able to show that $$\DT^{\omega_S}_h(v;1)=\s^{-p_g}\VW_h(v),$$ where $\VW_h(-)$ is the Vafa-Witten invariant (which was mathematically defined and studied in detail by Tanaka and Thomas in \cite{TT1,TT2}) and is expected to have modular properties based on S-duality conjecture (see \cite{VW94}).
% $\DT_h(v)$ by taking the virtual Euler number of the $\C^*$-fixed locus of $\mMb_h(v)$ (see Definiton \ref{DThv4} and Remarks \ref{virnum}, \ref{vafawitten}).

The $\C^*$-fixed locus $\mMw_h(v)^{\C^*}$ consists of sheaves supported on $S$ (the zero section of $\omega_S$) and its thickenings. One can write $\mMw_h(v)^{\C^*}$ as a disjoint union of several types of components, where each type is indexed by a partition of $r$. Out of these component types, there are two types of particular importance; One of them (we call it type I) is identified with $\mM_h(v)$, the moduli space of rank $r$ torsion free stable sheaves on $S$. The other type (we call it type II) can be identified with the nested Hilbert scheme $\S{\n}_\betab$ for a suitable choice of $\n, \betab$ depending on $v$. The reason that types I and II are more interesting, is the following result proven in \cite{GSY17b}:
\begin{theorem}(\cite[Theorem 2]{GSY17b})
The restriction of $[\mMw_h(v)^{\C^*}]^{\vir}_{\red}$ to the type I component $\mM_h(v)$ is identified with $[\mM_h(v)]^{\vir}_0$ induced by the natural trace free perfect obstruction theory over $\mM_h(v)$.
The restriction of $[\mMw_h(v)^{\C^*}]^{\vir}_{\red}$ to a type II component $\S{\n}_\betab$ is identified with  $[\S{\n}_\betab]^{\vir}$.
\end{theorem}
%Propositions \ref{(r)} and \ref{fixedpart}), the restriction of the fixed part of the perfect obstruction theory of $\mMb_h(v)$ to these component types respectively identified with the natural trace free perfect obstruction theory on $\mM_h(v)$, and the perfect obstruction theory we constructed for $\S{\n}_\betab$ in Theorem \ref{thm1}. Let $\chi^{\vir}(\mM_h(v))$ and $\chi^{\vir}(\S{\n}_{\betab})$ be the virtual Euler numbers with respect to these perfect obstruction theories \cite{FG10}. 

When $r=2$, then types I and II are the only component types of $\mMw_h(v)^{\C^*}$. This leads us to the following result: %(Proposition \ref{virdec}, \ref{(r)}, \ref{fixedpart}):

\begin{theorem}(\cite[Theorem 3]{GSY17b})
Suppose that $v=(2,\gamma,m)$. Then, 
\begin{align*} 
\DT^{\omega_S}_h(v;\alpha)&=\DT^{\omega_S}_h(v;\alpha)_{{\rm{I}}}+\sum_{n_1,n_2,\beta}\DT^{\omega_S}_h(v;\alpha)_{{\rm{II}}, \S{n_1,n_2}_\beta},\\
\DT^{\omega_S}_h(v)&=\chi^{\vir}(\mM_h(v))+\sum_{n_1,n_2,\beta}\chi^{\vir}(\S{n_1,n_2}_{\beta}),\end{align*} 
where the sum is over all $n_1, n_2, \beta$ (depending on $v$) for which $\S{n_1,n_2}_\beta$ is a type II component of  $\mMw_h(v)^{\C^*}$, and the indices I and II indicate the contributions of type I and II components to the invariant $\DT^{\omega_S}_h(v;\alpha)$.
\end{theorem}

The stability of sheaves imposes a strong condition on $n_1, n_2, \beta$  appearing in the summation in the theorem above. For example, if $S$ is a generic complete intersection in a projective space, then for any $n_1, n_2, \beta$ for which $\S{n_1,n_2}_\beta$ is a type II component of  $\mMw_h(v)^{\C^*}$, the condition in Theorem \ref{thm1.5} (leading to the vanishing $[\S{n_1,n_2}_\beta]^{\vir}=0$)  is not satisfied. 
 
 The invariants $\chi^{\vir}(\S{n_1,n_2}_{\beta})$ and $\DT^{\omega_S}_h(v;\alpha)_{{\rm{II}},\S{n_1,n_2}_\beta}$ (for a suitable choice of class $\alpha$ e.g. $\alpha=1$) appearing in the theorem above are special types of the invariants $$\sN_S(n_1,n_2,\beta;-)$$ which are defined as follows:
 
\begin{defi} \label{invs} ( \cite[Definition 2.13]{GSY17a}) Let $M\in \pic(S)$. Define the following elements in  $K(\S{n_1,n_2}_\beta)$ of virtual ranks respectively $n_1+n_2$ and $-\beta\cdot \beta^D/2+\beta\cdot c_1(M)$: $$\sK^{n_1,n_2}_{\beta;M}:=\left[\dR\pi_*M(\Z_\beta) \right]-\left[\dR\hom_\pi(\i{n_1},\i{n_2}_{\beta}\otimes M)\right], \quad \sG_{\beta;M}:=\left[\dR\pi_*M(\Z_\beta)|_{\Z_\beta}\right].$$  If $\beta=0$ we will instead use the notation $\sK^{[n_1\ge n_2]}_M:=\sK^{n_1,n_2}_{0;M}$ (see \cite[Definition 5.4]{GSY17a}). %We use the same notation for the restriction of $E_M$  to $\SD{\n} \subset \SD{n_1}\times \SD{n_2}$ in case $n_1\ge n_2$. Define the following generating series:
%We define the invariant
%$$\sN_S(n_1,n_2,\beta; M):=\int_{[\S{n_1,n_2}_\beta]^{\vir}}c_{n_1+n_2}(E^{n_1,n_2}_{\beta,M})\cup c_{e}(K_\beta).$$
We also define the rank $2n_i$ twisted tangent bundles $$\sT^M_{\S{n_i}}:=\left[\dR\pi_*M \right]-\left[\dR\hom_\pi(\i{n_i},\i{n_i}\otimes M)\right]=\left[\ext^1_\pi\left(\i{n_i},\i{n_i}\otimes M\right)_0\right].$$ Note that if $M=\O$ then $\sT^M_{\S{n_i}}=T_{\S{n_i}}$.

Let $\cP:=\cP(M,\beta,n_1,n_2)$ be a polynomial in the Chern classes of $\sK^{n_1,n_2}_{\beta;M}$, $\sG_{\beta;M}$, $\sT_{\S{n_1}}$, and $\sT_{\S{n_2}}$, then, we can define the invariant
$$\sN_S(n_1,n_2,\beta; \cP):=\int_{[\S{n_1,n_2}_\beta]^{\vir}}\cP.$$
Moreover under the condition in Theorem \ref{thm1.5}, the authors defined the reduced invariants $$\sN^{\red}_S(n_1,n_2,\beta; \cP):=\int_{[\S{n_1,n_2}_\beta]^{\vir}_{\red}}\cP.$$
\end{defi} 
 
Mochizuki in \cite{M02} expresses certain integrals against the virtual cycle of $\mM_h(v)$ in terms of Seiberg-Witten invariants and integrals $\sA(\gamma_1, \gamma_2, v;-)$ over the product of Hilbert scheme of points on $S$. Using this result the authors  were able to prove the following :

\begin{theorem} (\cite[Theorem 4]{GSY17b})
Suppose that  $p_g(S)>0$, and $v=(2,\gamma,m)$ is such that $\gamma\cdot h >2K_S \cdot h$ and $\chi(v) :=\int_S v \cdot td_S \ge 1$. Then,
\begin{align*}\label{Moc:2}
\DT^{\omega_S}_h(v;1)
=&
-\sum_{\begin{subarray}{c}
\gamma_1 + \gamma_2 =\gamma \\
\gamma_1\cdot h < \gamma_2 \cdot h
\end{subarray}}
\mathrm{SW}(\gamma_1) \cdot 2^{2-\chi(v)} \cdot \sA(\gamma_1, \gamma_2, v;\sP_1)+\sum_{n_1,n_2,\beta} \sN_S(n_1,n_2,\beta;\cP_1).\\
\DT^{\omega_S}_h(v)
=&
-\sum_{\begin{subarray}{c}
\gamma_1 + \gamma_2 =\gamma \\
\gamma_1\cdot h < \gamma_2 \cdot h
\end{subarray}}
\mathrm{SW}(\gamma_1) \cdot 2^{2-\chi(v)} \cdot \sA(\gamma_1, \gamma_2, v;\sP_2)+\sum_{n_1,n_2,\beta} \sN_S(n_1,n_2,\beta;\cP_2).
\end{align*}
Here $\mathrm{SW}(-)$ is the Seiberg-Witten invariant of $S$, $\sP_i$ and $\cP_i$ are certain universally defined explicit integrands, and the second sum in the formulas is over all $n_1, n_2, \beta$ (depending on $v$) for which $\S{n_1,n_2}_\beta$ is a type II component of  $\mMw_h(v)^{\C^*}$.

If $S$ is a $K3$ surface or $S$ is isomorphic to one of the five types of generic complete intersections $$(5)\subset \mathbb{P}^3, \; (3,3)\subset \mathbb{P}^4,\; (4,2)\subset \mathbb{P}^4, \; (3,2,2)\subset \mathbb{P}^5,\; (2,2,2,2)\subset \mathbb{P}^6,$$ the DT invariants $\DT^{\omega_S}_h(v;1)$ and $\DT^{\omega_S}_h(v)$ can be completely expressed as the sum of integrals over the product of the Hilbert schemes of points on $S$.

\end{theorem}

Center for Mathematical Sciences and Applications, Harvard University, Department of Mathematics,
20 Garden Street, Room 207, Cambridge, MA, 02139\\

Centre for Quantum Geometry of Moduli Spaces, Aarhus University, Department of Mathematics,
Ny Munkegade 118, building 1530, 319, 8000 Aarhus C, Denmark\\

National Research University Higher School of Economics, Russian Federation, Laboratory of Mirror Symmetry, NRU HSE, 6 Usacheva str., Moscow, Russia, 119048

\noindent {\tt{artan@cmsa.fas.harvard.edu}},

\end{document}